\documentclass{commat}

\title{%
On the Generalised Ricci Solitons and Sasakian Manifolds
    }
\author{%
A. Mohammed Cherif, K. Zegga and G. Beldjilali
}

\affiliation{
    \address{Ahmed Mohammed Cherif --
    Mascara University, Faculty of Exact Sciences, Department of Mathematics. Laboratory of Geometry, Analysis, Control and Applications. Mascara 29000, Algeria.
        }
    \email{%
    a.mohammedcherif@univ-mascara.dz
        }
    \address{Kaddour Zegga --
    Mascara University, Faculty of Exact Sciences, Department of Mathematics. Laboratory of Geometry, Analysis, Control and Applications. Mascara 29000, Algeria.
        }
    \email{%
    zegga.kadour@univ-mascara.dz
        }
    \address{Gherici Beldjlali --
    Mascara University, Faculty of Exact Sciences, Department of Mathematics. Laboratory of Quantum Physics and Mathematical Modeling (LPQ3M). Mascara 29000, Algeria.
        }
    \email{%
    gherici.beldjilali@univ-mascara.dz
        }
    }

\abstract{%
In this note, we find a necessary condition on odd-dimensional Riemannian manifolds under which both of Sasakian structure and the generalised Ricci soliton equation are satisfied, and we give some examples.
    }



\firstpage{119}

\VOLUME{30}

\DOI{https://doi.org/10.46298/cm.9311}

\begin{paper}

%%%%%%%%%%%%%%%%%%%%%%%%%%%%%%%%%%%%%%%
\section{Introduction and main results}
Let $(M,g)$ be a smooth Riemannian manifold.
By $R$ and $\operatorname{Ric}$
we denote respectively the Riemannian curvature tensor and the Ricci tensor of $(M,g)$.
Thus $R$ and $\operatorname{Ric}$ are defined by
\begin{equation}\label{eq1.1}
    R(X,Y)Z=\nabla_X \nabla_Y Z-\nabla_Y \nabla_X Z-\nabla_{[X,Y]}Z,
\end{equation}
\begin{equation}\label{eq1.2}
    \operatorname{Ric}(X,Y)=g(R(X,e_i)e_i,Y),
\end{equation}
where $\nabla$ is the Levi-Civita connection with respect to $g$, $\{e_i\}$ is an orthonormal frame, and $X,Y,Z\in \Gamma(TM)$.
 The gradient of a smooth function $f$ on $M$ is defined by
\begin{equation}\label{eq1.3}
    g(\operatorname{grad}f,X)=X(f),\quad \operatorname{grad}f=e_i(f)e_i,
\end{equation}
     where $X\in \Gamma(TM)$.
     The Hessian of $f$ is defined by
\begin{equation}\label{eq1.4}
    (\operatorname{Hess} f)(X,Y)=g(\nabla_X \operatorname{grad}f,Y),
\end{equation}
where $X,Y\in \Gamma(TM)$. For $X\in\Gamma(TM)$, we define $X^\flat\in\Gamma(T^*M)$ by
\begin{equation}\label{eq1.5}
X^\flat(Y)=g(X,Y).
\end{equation}
(For more details of previous definitions, see for example \cite{ON}).\\
The generalised Ricci soliton equation in Riemannian manifold $(M,g)$ is defined by (see \cite{NR})
\begin{equation}\label{eq1.6}
    \mathcal{L}_X g=-2c_1 X^\flat\odot X^\flat+2c_2\operatorname{Ric}+2\lambda g,
\end{equation}
where $X\in\Gamma(TM)$,
$\mathcal{L}_X g$ is the Lie-derivative of $g$ along $X$ given by
\begin{equation}\label{eq1.7}
   (\mathcal{L}_X g)(Y,Z)=g(\nabla_Y X,Z)+g(\nabla_Z X, Y),
\end{equation}
for all $Y,Z\in \Gamma(TM)$, and $c_1,c_2,\lambda\in \mathbb{R}$. Equation (\ref{eq1.6}), is a  generalization of
Killing's equation ($c_1=c_2=\lambda=0$),
 Equation for homotheties ($c_1=c_2=0$),
 Ricci soliton ($c_1=0$, $c_2=-1$),
 Cases of Einstein-Weyl ($c_1=1$, $c_2=\frac{-1}{n-2}$),
 Metric projective structures with
skew-symmetric Ricci tensor in projective class ($c_1=1$, $c_2=\frac{-1}{n-1}$, $\lambda=0$),
Vacuum near-horzion geometry equation ($c_1=1$, $c_2=\frac{1}{2}$), and is also a generalization of Einstein manifolds
(For more details, see  \cite{Besse}, \cite{Chrusciel}, \cite{Jacek}, \cite{Kunduri}, \cite{NR}).\\
 In this paper, we give a new generalization of Ricci soliton equation in Riemannian manifold $(M,g)$, given by the following equation
\begin{equation}\label{eq1.8}
    \mathcal{L}_{X_1} g=-2c_1 X_2^\flat\odot X_2^\flat+2c_2\operatorname{Ric}+2\lambda g,
\end{equation}
where $X_1,X_2\in\Gamma(TM)$.\\
Note that, if $X_1=\operatorname{grad}f_1$ and $X_2=\operatorname{grad}f_2$, where $f_1,f_2\in C^\infty(M)$, the generalised Ricci soliton equation (\ref{eq1.8})
is given by
\begin{equation}\label{eq1.9}
  \operatorname{Hess}f_1  =-c_1 df_2\odot df_2+c_2\operatorname{Ric}+\lambda g.
\end{equation}
\begin{example} Let $\mathbb{H}^2=\{(x,y)\in\mathbb{R}^2|y>0\}$ be a $2$-dimensional hyperbolic space equipped with the
Riemannian metric $g=\frac{dx^2+dy^2}{y^2}$, the following  functions $$f_1(x,y)=-(\lambda-c_2)\ln y,\quad f_2(x,y)=-\frac{\sqrt{c_1(\lambda-c_2)}}{c_1}\ln y,$$
 satisfy the generalised Ricci soliton equation (\ref{eq1.9}) with $c_1(\lambda-c_2)>0$.
\end{example}
\begin{example}\label{example2}
The product Riemannian manifold $M^3=(0,\infty)\times\mathbb{R}^2$
equipped with the Riemannian metric $g=dx^2+x^2(dy^2+dz^2)$ satisfies the generalised Ricci soliton equation (\ref{eq1.9}), with
$$f_1(x,y,z)=\frac{\lambda}{2} x^2-c_2\ln x,\quad f_2(x,y,z)=-\frac{\sqrt{-c_1c_2}}{c_1}\ln x,$$
where $c_1c_2<0$.
\end{example}
\begin{remark}
There are Riemannian manifolds that do not admit generalized soliton equation (\ref{eq1.9}) such that $f_1 = f_2$ (for example, the Riemannian manifold given in Example \ref{example2}).
\end{remark}

An $(2n+1)$-dimensional Riemannian manifold $(M,g)$ is said to be an almost contact
metric manifold if there exist on $M$ a $(1,1)$ tensor field $\varphi$, a vector
field $\xi$ (called the structure vector field) and a $1$-form $\eta$ such that
$$\eta(\xi)=1, \quad \varphi^{2}(X) = -X+\eta(X)\xi,  \quad g(\varphi X,\varphi Y) = g(X,Y)-\eta(X)\eta(Y),$$
 for any $X,Y\in\Gamma(TM)$. In particular, in an almost contact metric manifold we also
have $\varphi\xi=0$ and $\eta \circ \varphi=0$.
Such a manifold is said to be a contact metric manifold if $d\eta=\phi$, where
$\phi(X,Y)=g(X,\varphi Y)$ is called the fundamental $2$-form of $M$. If, in addition,
$\xi$ is a Killing vector field, then $M$ is said to be a K-contact manifold. It is
well-known that a contact metric manifold is a K-contact manifold if and only if
$ \nabla_{X}\xi = -\varphi X $, for any vector field $X$ on $M$.
The almost contact metric structure of $M$ is said to be normal if
 $[\varphi,\varphi](X,Y)=-2d\eta~(X,Y)\xi$, for any $X,Y\in\Gamma(TM)$, where $[\varphi,\varphi]$ denotes the
Nijenhuis torsion of $\varphi$, given by
$$ [\varphi,\varphi](X,Y)=\varphi^{2}[X,Y]+[\varphi X,\varphi Y]-\varphi[\varphi X,Y]-\varphi[X,\varphi Y]. $$
A normal contact metric manifold is called a Sasakian manifold. It can be
proved that a Sasakian manifold is K-contact, and that an almost contact metric
manifold is Sasakian if and only if
\begin{equation}\label{FormSas1}
(\nabla_{X}\varphi)Y = g(X,Y)\xi-\eta(Y)X,
\end{equation}
for any $X$,$Y$. Moreover, for a Sasakian manifold the following equation holds
$$ R(X,Y)\xi= \eta(Y)X - \eta(X)Y .$$
From the formula (\ref{FormSas1}) easily obtains
\begin{equation}\label{FormSas2}
  \nabla_{X}\xi = - \varphi X, \qquad   (\nabla_{X}\eta)Y = -g(\varphi X,Y).
\end{equation}
(For more details, see \cite{BL1}, \cite{BG}, \cite{YK}).\\
The main result of this paper is the following:
\begin{theorem}\label{thm}
Suppose $(M,\varphi , \xi , \eta, g)$ is a Sasakian manifold, and satisfies the generalised Ricci soliton
equation $(\ref{eq1.9})$. Then
\begin{equation}\label{condition}
\zeta\equiv\operatorname{grad}f_1+c_1\xi(\xi(f_2))\operatorname{grad}f_2-c_1\xi(f_2)\nabla_\xi \operatorname{grad}f_2=\xi(f_1)\xi.
\end{equation}
\end{theorem}
\begin{remark}
The condition (\ref{condition}) is necessary for the existence of a Sasakian structure and the generalised Ricci soliton
equation $(\ref{eq1.9})$ on an odd-dimensional Riemannian manifold.
\end{remark}
\begin{example}
Consider the Sasakian manifold $(\mathbb{R}^2\times(0,\pi),\varphi , \xi , \eta, g)$ endowed with the Sasakian structure $(\varphi , \xi , \eta, g)$ given by
 $$(g_{ij})=\left(
                                                                                                       \begin{array}{ccc}
                                                                                                         p^2+q^2 & 0 & -q \\
                                                                                                         0 & p^2 & 0 \\
                                                                                                         -q & 0 & 1 \\
                                                                                                       \end{array}
                                                                                                     \right),\quad
 (\varphi_{ij})=\left(
                                                                                                       \begin{array}{ccc}
                                                                                                         0 & -1 & 0 \\
                                                                                                         1 & 0 & 0 \\
                                                                                                         0 & -q & 0 \\
                                                                                                       \end{array}
                                                                                                     \right),$$
$$\xi=\frac{\partial}{\partial z},\quad\eta=-qdx+dz,\quad p(x,y,z)=\frac{4 e^y}{16+e^{2y}},\quad q(x,y,z)=\frac{- e^{2y}}{16+e^{2y}}.$$
Then, the following smooth functions
$$f_1(x,y,z)=\frac{2c_2+\lambda}{2}\left(\ln(16+e^{2y})-2\ln\left(\frac{\sin z}{2c_2+\lambda}\right)\right),$$
$$f_2(x,y,z)=-\frac{1}{2}\sqrt{-\frac{2c_2+\lambda}{c_1}}\left(2\ln(\sin z)-\ln(16+e^{2y})\right),$$
satisfy the generalised Ricci soliton equation (\ref{eq1.9}), where $c_1<0$ and $2c_2+\lambda>0$. Furthermore,
$$\zeta=\xi(f_1)\xi=-(2c_2+\lambda)\cot(z)\xi.$$
\end{example}

%%%%%%%%%%%%%%%%%%%%%%%%%%%%%%%%%%%%%%%%%%%%%%%%%%%%%%%%%%%%%%%%%%%%%%%%%%

\section{Proof of the result}
%
%
%
%
%
%
%
%
%
%
%
%%%%%%%%%%%%%%%%%%%%%%%%%%%%%%%%%%%%%%%%%%%%%%%%%%%%%%%%%%%%%%%%%%%%%%%%%%%%%%%%%%
%%%%%%%%%%%%%%%%%%%%%%%%%% Contents of Section 2 %%%%%%%%%%%%%%%%%%%%%%%%%%%%%%%%%
%%%%%%%%%%%%%%%%%%%%%%%%%%%%%%%%%%%%%%%%%%%%%%%%%%%%%%%%%%%%%%%%%%%%%%%%%%%%%%%%%%
For the proof of Theorem \ref{thm}, we need the following lemmas.

\begin{lemma}\label{lem1}\cite{mc}
Let $(M,\varphi , \xi , \eta, g)$ be a  Sasakian manifold. Then
\begin{eqnarray*}
% \nonumber to remove numbering (before each equation)
  \big(\mathcal{L}_\xi(\mathcal{L}_{X_1}g)\big)(Y,\xi)
   &=& g(X_1,Y)+g(\nabla_\xi\nabla_\xi X_1,Y) +Y g(\nabla_\xi X_1, \xi),
\end{eqnarray*}
where $X_1,Y\in\Gamma(TM)$, with $Y$ is orthogonal to $\xi$.
\end{lemma}
\begin{lemma}\label{lem2}\cite{mc}
Let $(M,g)$ be a Riemannian manifold, and let $f_2\in C^\infty(M)$. Then
\begin{eqnarray*}
% \nonumber to remove numbering (before each equation)
  \big(\mathcal{L}_\xi(df_2\odot df_2)\big)(Y,\xi)
   &=&Y(\xi(f_2))\xi(f_2)+Y(f_2)\xi(\xi(f_2)),
\end{eqnarray*}
where $\xi,Y\in\Gamma(TM)$.
\end{lemma}
\begin{lemma}\label{lem3}
Let $(M,\varphi , \xi , \eta, g)$ be a Sasakian manifold of dimension $(2n + 1)$, and
satisfies the generalised Ricci soliton equation $(\ref{eq1.9})$. Then
\begin{eqnarray*}
\nabla_\xi\operatorname{grad}f_1=(\lambda+2c_2n)\xi-c_1\xi(f_2)\operatorname{grad}f_2.
\end{eqnarray*}
\end{lemma}
\begin{proof} Let $Y\in\Gamma(TM)$, we have
\begin{eqnarray*}
\operatorname{Ric}(\xi,Y)
&=&g(R(\xi,e_i)e_i,Y)\\
&=&g(R(e_i,Y)\xi,e_i)\\
&=&\eta(Y)g(e_i,e_i)-\eta(e_i)g(X,e_i)\\
&=&(2n+1)\eta(Y)-\eta(Y)\\
&=&2n\eta(Y)\\
&=&2ng(\xi,Y),
\end{eqnarray*}
where $\{e_i\}$ is an orthonormal frame on $M$, which implies
\begin{eqnarray}\label{eq2.5}
% \nonumber to remove numbering (before each equation)
\lambda g(\xi,Y)+c_2\operatorname{Ric}(\xi,Y)
   &=&\nonumber \lambda g(\xi,Y)+2c_2ng(\xi,Y) \\
   &=& (\lambda+2c_2n)g(\xi,Y).
\end{eqnarray}
From equations (\ref{eq1.9}) and (\ref{eq2.5}), we obtain
\begin{eqnarray}\label{eq2.6}
% \nonumber to remove numbering (before each equation)
(\operatorname{Hess}f_1)(\xi,Y)
  &=&\nonumber  -c_1\xi(f_2)Y(f_2) +(\lambda+2c_2n)g(\xi,Y)\\
   &=&  -c_1\xi(f_2)g(\operatorname{grad}f_2,Y) +(\lambda+2c_2n)g(\xi,Y),
\end{eqnarray}
the Lemma follows from equation (\ref{eq2.6}).
\end{proof}
\begin{proof}
[Proof of Theorem \ref{thm}] Let $Y\in\Gamma(TM)$, such that $g(\xi,Y)=0$, from Lemma \ref{lem1}, with $X_1= \operatorname{grad}f_1$, we have
\begin{eqnarray}\label{eq2.7}
% \nonumber to remove numbering (before each equation)
2\big(\mathcal{L}_\xi(\operatorname{Hess}f_1)\big)(Y,\xi)
   &=&\nonumber Y(f_1)+g(\nabla_\xi \nabla_\xi \operatorname{grad}f_1,Y)\\
   &&+Yg(\nabla_\xi  \operatorname{grad}f_1,\xi).
\end{eqnarray}
By Lemma \ref{lem3}, and equation (\ref{eq2.7}), we get
\begin{eqnarray}\label{eq2.8}
% \nonumber to remove numbering (before each equation)
2\big(\mathcal{L}_\xi(\operatorname{Hess}f_1)\big)(Y,\xi)
   &=&\nonumber Y(f_1)+(\lambda+2c_2n)g(\nabla_\xi\xi,Y)\\
   &&\nonumber-c_1g(\nabla_\xi (\xi(f_2)\operatorname{grad}f_2),Y)  \\
   & &  +(\lambda+2c_2n)Yg(\xi,\xi)-c_1Y(\xi(f_2)^2).
\end{eqnarray}
Since $\nabla_\xi\xi=0$ and $g(\xi,\xi)=1$, from equation (\ref{eq2.8}), we obtain
\begin{eqnarray}\label{eq2.9}
% \nonumber to remove numbering (before each equation)
2\big(\mathcal{L}_\xi(\operatorname{Hess}f_1)\big)(Y,\xi)
   &=&\nonumber Y(f_1)-c_1\xi(\xi(f_2))Y(f_2)
   \\
   &&\nonumber-c_1\xi(f_2)g(\nabla_\xi \operatorname{grad}f_2,Y)\\
   &&-2c_1\xi(f_2)Y(\xi(f_2)).
\end{eqnarray}
Since $\mathcal{L}_\xi g=0$ (i.e. $\xi$ is a Killing vector field), it implies that $\mathcal{L}_\xi \operatorname{Ric}=0$. Taking the Lie derivative
to the generalised Ricci soliton equation (\ref{eq1.9}) yields
\begin{equation}\label{eq2.11}
   2\big(\mathcal{L}_\xi(\operatorname{Hess}f_1)\big)(Y,\xi)
   =-2c_1 \big(\mathcal{L}_\xi(df_2\odot df_2)\big)(Y,\xi).
\end{equation}
Thus, from equations (\ref{eq2.9}), (\ref{eq2.11}) and Lemma \ref{lem2}, we have
\begin{eqnarray}\label{eq2.12}
% \nonumber to remove numbering (before each equation)
 &&\nonumber Y(f_1)-c_1\xi(\xi(f_2))Y(f_2)
   \\
   &&\nonumber-c_1\xi(f_2)g(\nabla_\xi \operatorname{grad}f_2,Y)-2c_1\xi(f_2)Y(\xi(f_2))\\
   &&=-2c_1 Y(\xi(f_2))\xi(f_2)-2c_1Y(f_2)\xi(\xi(f_2)),
\end{eqnarray}
which is equivalent to
\begin{equation}\label{eq2.13}
    Y(f_1)+c_1\xi(\xi(f_2))Y(f_2) -c_1\xi(f_2)g(\nabla_\xi \operatorname{grad}f_2,Y)=0,
\end{equation}
that is, the vector field
\begin{equation}\label{eq2.14}
\zeta=\operatorname{grad}f_1+c_1\xi(\xi(f_2))\operatorname{grad}f_2-c_1\xi(f_2)\nabla_\xi \operatorname{grad}f_2,
\end{equation}
is parallel to $\xi$. The proof is completed.
\end{proof}

%%% References

\EditInfo{%
    July 27, 2019}{%
    June 30, 2020}{%
    Haizhong Li}

\end{paper}